\documentclass[11pt]{article}

\usepackage{amsmath,amsfonts,amsthm,amssymb,graphicx,tikz,tikz-cd,url,hyperref,cleveref}
\usepackage[all,2cell,ps]{xy}

\theoremstyle{plain}
\newtheorem{thm}{Theorem}[section]

\newtheorem{lem}[thm]{Lemma}

\newtheorem{cor}[thm]{Corollary}

\newtheorem{qtn}[thm]{Question}

\theoremstyle{definition}

\newtheorem{rem}[thm]{Remark}

\theoremstyle{remark}

\newcommand{\bbH}{\mathbb{H}}

\newcommand{\bbP}{\mathbb{P}}

\newcommand{\bbR}{\mathbb{R}}
\newcommand{\bbS}{\mathbb{S}}

\newcommand{\bbZ}{\mathbb{Z}}

\newcommand{\calS}{\mathcal{S}}

\newcommand{\Del}{\Delta}

\newcommand{\thet}{\theta}

\DeclareMathOperator{\Orth}{O}

\newenvironment{pf}{\begin{proof}}{\end{proof}}

\usepackage{tikz-cd}
\allowdisplaybreaks

\makeatletter
\let\@@pmod\pmod
\DeclareRobustCommand{\pmod}{\@ifstar\@pmods\@@pmod}
\def\@pmods#1{\mkern4mu({\operator@font mod}\mkern 6mu#1)}
\makeatother

\usetikzlibrary{arrows}

\title{A hyperbolic $4$-orbifold with underlying space $\bbP^2$}
\author{Matthew Stover \\ \small{Temple University}\\ \small{\textsf{mstover@temple.edu}}}
\date{\today}

\begin{document}

\maketitle

\begin{abstract}
This paper shows that the complex projective plane $\bbP^2$ can be realized as the underlying space for a closed hyperbolic $4$-orbifold. This is the first example of a closed hyperbolic $4$-orbifold whose underlying space is symplectic, which is related to the open question as to whether or not closed hyperbolic $4$-manifolds can admit symplectic structures.
\end{abstract}

\section{Introduction}\label{sec:Intro}

One of the major open problems in higher-dimensional hyperbolic geometry is whether or not there is a closed hyperbolic $4$-manifold that admits a symplectic structure. LeBrun famously conjectured that the Seiberg--Witten invariants of closed hyperbolic $4$-manifolds vanish \cite{LeBrun}, which would immediately imply that closed hyperbolic $4$-manifolds cannot be symplectic. Lin and Martelli recently proved this conjecture for the Davis manifold \cite{LinMartelli}; see that paper for additional related references. This note shows that the orbifold version of LeBrun's conjecture is false.

\begin{thm}\label{thm:Main}
There is a complete hyperbolic $4$-orbifold with underlying space the complex projective plane $\bbP^2$. In particular, there is a closed symplectic hyperbolic $4$-orbifold.
\end{thm}

The proof, given in \Cref{sec:Construction}, is giving by producing a simplicial decomposition of $\bbP^2$ that can be realized by $60$ copies of a particular compact Coxeter simplex in $\bbH^4$. The orbifold locus is described in \Cref{sec:Locus}. I also found examples with underlying space the product $\bbS^2 \times \bbS^2$ of two spheres by similar methods, and Bruno Martelli described a convincing heuristic to me that one should also be able to construct a hyperbolic $4$-orbifold with underlying space the $4$-torus. The reason for focusing on $\bbP^2$ for this note is the following consequence, which is perhaps surprising in comparison with the fact that closed hyperbolic $4$-manifolds have signature zero by the Hirzebruch signature theorem \cite[\S 19]{Milnor}.

\begin{cor}\label{cor:SigCor}
There is a compact complete hyperbolic $4$-orbifold whose underlying space is a $4$-manifold with nonzero signature.
\end{cor}

The $4$-sphere is well-known to experts to be the underlying space of a hyperbolic $4$-orbifold by taking the orientation-preserving double cover of a right-angled reflection orbifold. Having realized arguably the simplest simply connected $4$-manifolds, $\bbS^4$, $\bbP^2$, and $\bbS^2 \times \bbS^2$, the following question seems quite interesting.

\begin{qtn}
Let $N$ be a simply connected closed $4$-manifold. Is $N$ the underlying space for a complete hyperbolic $4$-orbifold?
\end{qtn}

Lastly, the Coxeter group in $\Orth^+(4,1)$ generated by reflections in the sides of the simplex $\calS$ in \Cref{sec:Construction} is arithmetic \cite{Lanner, Vinberg}. Thus the following corollary to \Cref{thm:Main} is immediate from the method of construction.

\begin{cor}\label{cor:Arithmetic}
The orbifold constructed in this paper to prove \Cref{thm:Main} is arithmetic.
\end{cor}

\subsubsection*{Acknowledgments}
I thank for Bruno Martelli for conversations related to this paper and his observation that this example was particularly interesting because the signature of the underlying space is nonzero. I also thank Jonathan Spreer for comments on smooth structures. The author acknowledges support from the Institut Henri Poincar\'e (UAR 839 CNRS-Sorbonne Universit\'e) and LabEx CARMIN (ANR-10-LABx-59-01). This material is based upon work supported by NSF grants DMS-2203555 and DMS-2506896, along with award SFI-MPS-TSM-00014184 from the Simons Foundation. Part of this research was performed while the author was visiting the Mathematical Sciences Research Institute (MSRI), now becoming the Simons Laufer Mathematical Sciences Institute (SLMath), which is supported by the National Science Foundation (Grant No.\ DMS-1928930).

\section{A simplicial decomposition of $\bbP^2$}\label{sec:Construction}

This section describes a simplicial decomposition of $\bbP^2$ with $60$ simplices that equips it with a complete hyperbolic orbifold metric. Consider the oriented standard $4$-simplex $\Del$ in $\bbR^4$ with basis elements $x_0, \dots, x_4$. The notation $\Del(a_0\ \cdots\ a_d)$ will denote the $d$-dimensional facet of $\Del$ defined by the basis elements $a_j \in \{x_0, \dots, x_4\}$ with orientation associated with the given ordering. Number sixty $4$-simplices as $\Del_0, \dots, \Del_{59}$.

This simplicial decomposition will have the property that if $\Del_j$ is glued to $\Del_k$ along the face $\Del_j(a_0 \cdots a_3)$, then the corresponding face of $\Del_k$ is $\Del_k(a_0 \cdots a_3)$. Thus Tables \ref{tb:P2Decomp} and \ref{tb:P2Decomp2} define a unique simplicial complex, where $\Del_k$ in the row labeled $\Del_j$ and column labeled $(a_0 \cdots a_3)$ means that $\Del_j(a_0 \cdots a_3)$ is glued to $\Del_k(a_0 \cdots a_3)$; note that each gluing is therefore by a reflection through the given face.

\begin{table}[h!]
\centering
\begin{tabular}{|c||c|c|c|c|c|}
\hline
& $(0\, 1\, 2\, 3)$ & $(0\, 1\, 2\, 4)$ & $(0\, 1\, 3\, 4)$ & $(0\, 2\, 3\, 4)$ & $(1\, 2\, 3\, 4)$ \\
\hline\hline
$\Del_0$ & $\Del_1$ & $\Del_1$ & $\Del_2$ & $\Del_3$ & $\Del_2$ \\
$\Del_1$ & $\Del_0$ & $\Del_0$ & $\Del_4$ & $\Del_5$ & $\Del_6$ \\
$\Del_2$ & $\Del_4$ & $\Del_6$ & $\Del_0$ & $\Del_7$ & $\Del_0$ \\
$\Del_3$ & $\Del_5$ & $\Del_5$ & $\Del_8$ & $\Del_0$ & $\Del_9$ \\
$\Del_4$ & $\Del_2$ & $\Del_{10}$ & $\Del_1$ & $\Del_{11}$ & $\Del_{10}$ \\
$\Del_5$ & $\Del_3$ & $\Del_3$ & $\Del_{12}$ & $\Del_1$ & $\Del_{13}$ \\
$\Del_6$ & $\Del_{10}$ & $\Del_2$ & $\Del_{10}$ & $\Del_{14}$ & $\Del_1$ \\
$\Del_7$ & $\Del_{11}$ & $\Del_{14}$ & $\Del_{15}$ & $\Del_2$ & $\Del_{16}$ \\
$\Del_8$ & $\Del_{12}$ & $\Del_{17}$ & $\Del_3$ & $\Del_{15}$ & $\Del_{18}$ \\
$\Del_9$ & $\Del_{19}$ & $\Del_{13}$ & $\Del_{18}$ & $\Del_{16}$ & $\Del_3$ \\
$\Del_{10}$ & $\Del_6$ & $\Del_4$ & $\Del_6$ & $\Del_{20}$ & $\Del_4$ \\
$\Del_{11}$ & $\Del_7$ & $\Del_{20}$ & $\Del_{21}$ & $\Del_4$ & $\Del_{22}$ \\
$\Del_{12}$ & $\Del_8$ & $\Del_{23}$ & $\Del_5$ & $\Del_{21}$ & $\Del_{24}$ \\
$\Del_{13}$ & $\Del_{23}$ & $\Del_9$ & $\Del_{24}$ & $\Del_{25}$ & $\Del_5$ \\
$\Del_{14}$ & $\Del_{20}$ & $\Del_7$ & $\Del_{26}$ & $\Del_6$ & $\Del_{25}$ \\
$\Del_{15}$ & $\Del_{21}$ & $\Del_{27}$ & $\Del_7$ & $\Del_8$ & $\Del_{28}$ \\
$\Del_{16}$ & $\Del_{29}$ & $\Del_{25}$ & $\Del_{28}$ & $\Del_9$ & $\Del_7$ \\
$\Del_{17}$ & $\Del_{24}$ & $\Del_8$ & $\Del_{23}$ & $\Del_{27}$ & $\Del_{30}$ \\
$\Del_{18}$ & $\Del_{30}$ & $\Del_{30}$ & $\Del_9$ & $\Del_{31}$ & $\Del_8$ \\
$\Del_{19}$ & $\Del_9$ & $\Del_{24}$ & $\Del_{30}$ & $\Del_{29}$ & $\Del_{23}$ \\
$\Del_{20}$ & $\Del_{14}$ & $\Del_{11}$ & $\Del_{32}$ & $\Del_{10}$ & $\Del_{33}$ \\
$\Del_{21}$ & $\Del_{15}$ & $\Del_{34}$ & $\Del_{11}$ & $\Del_{12}$ & $\Del_{35}$ \\
$\Del_{22}$ & $\Del_{36}$ & $\Del_{33}$ & $\Del_{35}$ & $\Del_{33}$ & $\Del_{11}$ \\
$\Del_{23}$ & $\Del_{13}$ & $\Del_{12}$ & $\Del_{17}$ & $\Del_{34}$ & $\Del_{19}$ \\
$\Del_{24}$ & $\Del_{17}$ & $\Del_{19}$ & $\Del_{13}$ & $\Del_{37}$ & $\Del_{12}$ \\
$\Del_{25}$ & $\Del_{34}$ & $\Del_{16}$ & $\Del_{38}$ & $\Del_{13}$ & $\Del_{14}$ \\
$\Del_{26}$ & $\Del_{32}$ & $\Del_{39}$ & $\Del_{14}$ & $\Del_{32}$ & $\Del_{38}$ \\
$\Del_{27}$ & $\Del_{37}$ & $\Del_{15}$ & $\Del_{39}$ & $\Del_{17}$ & $\Del_{40}$ \\
$\Del_{28}$ & $\Del_{41}$ & $\Del_{40}$ & $\Del_{16}$ & $\Del_{42}$ & $\Del_{15}$ \\
$\Del_{29}$ & $\Del_{16}$ & $\Del_{37}$ & $\Del_{41}$ & $\Del_{19}$ & $\Del_{36}$ \\
\hline
\end{tabular}
\caption{The decomposition of $\bbP^2$ (ctd.\ in \Cref{tb:P2Decomp2})}\label{tb:P2Decomp}
\end{table}
\begin{table}[h!]
\centering
\begin{tabular}{|c||c|c|c|c|c|}
\hline
& $(0\, 1\, 2\, 3)$ & $(0\, 1\, 2\, 4)$ & $(0\, 1\, 3\, 4)$ & $(0\, 2\, 3\, 4)$ & $(1\, 2\, 3\, 4)$ \\
\hline\hline
$\Del_{30}$ & $\Del_{18}$ & $\Del_{18}$ & $\Del_{19}$ & $\Del_{43}$ & $\Del_{17}$ \\
$\Del_{31}$ & $\Del_{43}$ & $\Del_{43}$ & $\Del_{42}$ & $\Del_{18}$ & $\Del_{42}$ \\
$\Del_{32}$ & $\Del_{26}$ & $\Del_{44}$ & $\Del_{20}$ & $\Del_{26}$ & $\Del_{45}$ \\
$\Del_{33}$ & $\Del_{46}$ & $\Del_{22}$ & $\Del_{45}$ & $\Del_{22}$ & $\Del_{20}$ \\
$\Del_{34}$ & $\Del_{25}$ & $\Del_{21}$ & $\Del_{44}$ & $\Del_{23}$ & $\Del_{46}$ \\
$\Del_{35}$ & $\Del_{47}$ & $\Del_{46}$ & $\Del_{22}$ & $\Del_{48}$ & $\Del_{21}$ \\
$\Del_{36}$ & $\Del_{22}$ & $\Del_{48}$ & $\Del_{47}$ & $\Del_{46}$ & $\Del_{29}$ \\
$\Del_{37}$ & $\Del_{27}$ & $\Del_{29}$ & $\Del_{49}$ & $\Del_{24}$ & $\Del_{48}$ \\
$\Del_{38}$ & $\Del_{44}$ & $\Del_{50}$ & $\Del_{25}$ & $\Del_{49}$ & $\Del_{26}$ \\
$\Del_{39}$ & $\Del_{49}$ & $\Del_{26}$ & $\Del_{27}$ & $\Del_{44}$ & $\Del_{50}$ \\
$\Del_{40}$ & $\Del_{51}$ & $\Del_{28}$ & $\Del_{50}$ & $\Del_{52}$ & $\Del_{27}$ \\
$\Del_{41}$ & $\Del_{28}$ & $\Del_{51}$ & $\Del_{29}$ & $\Del_{53}$ & $\Del_{47}$ \\
$\Del_{42}$ & $\Del_{53}$ & $\Del_{52}$ & $\Del_{31}$ & $\Del_{28}$ & $\Del_{31}$ \\
$\Del_{43}$ & $\Del_{31}$ & $\Del_{31}$ & $\Del_{53}$ & $\Del_{30}$ & $\Del_{52}$ \\
$\Del_{44}$ & $\Del_{38}$ & $\Del_{32}$ & $\Del_{34}$ & $\Del_{39}$ & $\Del_{54}$ \\
$\Del_{45}$ & $\Del_{54}$ & $\Del_{54}$ & $\Del_{33}$ & $\Del_{55}$ & $\Del_{32}$ \\
$\Del_{46}$ & $\Del_{33}$ & $\Del_{35}$ & $\Del_{54}$ & $\Del_{36}$ & $\Del_{34}$ \\
$\Del_{47}$ & $\Del_{35}$ & $\Del_{56}$ & $\Del_{36}$ & $\Del_{56}$ & $\Del_{41}$ \\
$\Del_{48}$ & $\Del_{56}$ & $\Del_{36}$ & $\Del_{55}$ & $\Del_{35}$ & $\Del_{37}$ \\
$\Del_{49}$ & $\Del_{39}$ & $\Del_{57}$ & $\Del_{37}$ & $\Del_{38}$ & $\Del_{55}$ \\
$\Del_{50}$ & $\Del_{57}$ & $\Del_{38}$ & $\Del_{40}$ & $\Del_{57}$ & $\Del_{39}$ \\
$\Del_{51}$ & $\Del_{40}$ & $\Del_{41}$ & $\Del_{57}$ & $\Del_{58}$ & $\Del_{56}$ \\
$\Del_{52}$ & $\Del_{58}$ & $\Del_{42}$ & $\Del_{58}$ & $\Del_{40}$ & $\Del_{43}$ \\
$\Del_{53}$ & $\Del_{42}$ & $\Del_{58}$ & $\Del_{43}$ & $\Del_{41}$ & $\Del_{58}$ \\
$\Del_{54}$ & $\Del_{45}$ & $\Del_{45}$ & $\Del_{46}$ & $\Del_{59}$ & $\Del_{44}$ \\
$\Del_{55}$ & $\Del_{59}$ & $\Del_{59}$ & $\Del_{48}$ & $\Del_{45}$ & $\Del_{49}$ \\
$\Del_{56}$ & $\Del_{48}$ & $\Del_{47}$ & $\Del_{59}$ & $\Del_{47}$ & $\Del_{51}$ \\
$\Del_{57}$ & $\Del_{50}$ & $\Del_{49}$ & $\Del_{51}$ & $\Del_{50}$ & $\Del_{59}$ \\
$\Del_{58}$ & $\Del_{52}$ & $\Del_{53}$ & $\Del_{52}$ & $\Del_{51}$ & $\Del_{53}$ \\
$\Del_{59}$ & $\Del_{55}$ & $\Del_{55}$ & $\Del_{56}$ & $\Del_{54}$ & $\Del_{57}$ \\
\hline
\end{tabular}
\caption{The decomposition of $\bbP^2$ (ctd.)}\label{tb:P2Decomp2}
\end{table}

\begin{thm}\label{thm:CP2}
The simplicial complex described in Tables \ref{tb:P2Decomp} and \ref{tb:P2Decomp2} is a PL manifold that is smoothable to the standard smooth structure on $\bbP^2$.
\end{thm}

\begin{pf}
The computational software Regina \cite{Regina} can take a simplicial decomposition described as in this section and compute the homology groups and fundamental group. Regina quickly verifies that the given simplicial decomposition gives a simply connected manifold with $H_2(M) \cong \bbZ$, which must be $\bbP^2$ by Freedman's theorem \cite{Freedman}. In fact, several iterations of Regina's simplification algorithm quickly simplify the given simplicial decomposition to the four-simplex simplicial decomposition of $\bbP^2$ contained in Regina's database of example manifolds. As described by Burke, Burton, and Spreer \cite[\S 4.4]{BBS}, this simplicial decomposition can be obtained by applying Regina's standard simplifications, which do not change the underlying smooth structure, from a standard Kirby diagram for $\bbP^2$ associated with the standard smooth structure. This proves the theorem.
\end{pf}

\begin{rem}
A github repository for this paper \cite{StoverCode} contains a Regina file the reader can use to replicate the computations used to prove \Cref{thm:CP2}. The repository also contains a file with the matrices defining the induced chain complex, from which one can directly compute the homology groups using any linear algebra software.
\end{rem}

Now consider the compact simplex $\calS$ in $\bbH^4$ with Coxeter diagram
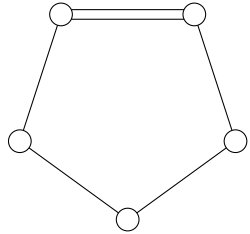
\[
\begin{tikzpicture}[scale=0.5]
\draw (-1.76336, 2.42705) -- (-2.85317, -0.927051) -- (0, -3) -- (2.85317, -0.927051) -- (1.76336, 2.42705);
\draw (-1.76336, 2.3) -- (1.76336, 2.3);
\draw (-1.76336, 2.5541) -- (1.76336, 2.5541);
\draw[fill=white] (0,-3) circle (0.3cm);
\draw[fill=white] (-1.76336, 2.42705) circle (0.3cm);
\draw[fill=white] (-2.85317, -0.927051) circle (0.3cm);
\draw[fill=white] (2.85317, -0.927051) circle (0.3cm);
\draw[fill=white] (1.76336, 2.42705) circle (0.3cm);
\end{tikzpicture}
\]
first found by Lann\'er \cite{Lanner}. The next goal is to prove the following more precise version of \Cref{thm:Main}.

\begin{thm}\label{thm:Orbifold} The given simplicial decomposition of $\bbP^2$ in fact is an orbifold tiling by $\calS$. In other words, appropriately identifying each $\Del_j$ with $\calS$ equips $\bbP^2$ with a complete hyperbolic orbifold metric.
\end{thm}

\begin{pf}
By the Poincar\'e polyhedron theorem \cite[Thm.\ 13.5.3]{Ratcliffe}, it suffices to show that the number of simplices around a vertex of type $(a_0\ a_1\ a_2)$ divides $2 \pi / \thet(a_0\ a_1\ a_2)$, where $\thet(a_0\ a_1\ a_2)$ is the dihedral angle of $\calS$ at that triangle. Specifically,
\begin{align*}
\theta(0\ 1\ 2) &= \frac{\pi}{4} & \theta(0\ 1\ 4) &= \frac{\pi}{3} & \theta(0\ 3\ 4) &= \frac{\pi}{3} \\
\theta(2\ 3\ 4) &= \frac{\pi}{3} & \theta(1\ 2\ 3) &= \frac{\pi}{3} & &
\end{align*}
and all other dihedral angles are $\pi / 2$. An accounting of the the triangles with label $(0\ 1\ 2)$ is in \Cref{tb:012}, and every degree is divisible by $8$, as required; interior points of a triangle with degree $8$ have neighborhoods isometric to $\bbH^4$, and if the degree is $d < 8$, then any interior point is contained in the orbifold locus with weight $8 / d$. Every triangle in the simplicial decomposition with the same labeling as a triangle of $\calS$ with dihedral angle $\pi / 3$ has exactly six simplices that meet at the given triangle, which makes the identification space locally isometric to $\bbH^4$ around points in the interior of those triangles. For the remaining triangles, where $\calS$ has dihedral angle $\pi / 2$, there are exactly four simplices around each triangle except those enumerated in \Cref{tb:BadTri}, which have only two; those with four simplices have neighborhoods isometric to $\bbH^4$ and those with two are part of the orbifold locus and have weight two. This completes all the checks necessary to apply the Poincar\'e polyhedron theorem.
\end{pf}

\begin{table}[h!]
\centering
\begin{tabular}{|c|c|c|c|c|}
\hline
Triangle & Identified $\Del_j$ & Degree & Edges & Vertices \\
\hline
$t_9$ & \begin{tabular}{c} $\Del_0$ $\Del_1$\end{tabular} & $2$ & $e_0, e_1, e_4$ & $v_0, v_1, v_2$ \\
\hline
$t_{15}$ & \begin{tabular}{c} $\Del_2$ $\Del_4$ $\Del_{10}$ $\Del_6$ \end{tabular} & $4$ & $e_0, e_{11}, e_4$ & $v_0, v_1, v_2$ \\
\hline
$t_{21}$ & \begin{tabular}{c} $\Del_3$ $\Del_5$ \end{tabular} & $2$ & $e_{12}, e_1, e_{13}$ & $v_0, v_5, v_2$ \\
\hline
$t_{34}$ & \begin{tabular}{c} $\Del_7$ $\Del_{14}$ $\Del_{20}$ $\Del_{11}$ \end{tabular} & $4$ & $e_{18}, e_{11}, e_{19}$ & $v_0, v_6, v_2$ \\
\hline
$t_{40}$ & \begin{tabular}{c} $\Del_8$ $\Del_{17}$ $\Del_{24}$ $\Del_{19}$ \\ $\Del_9$ $\Del_{13}$ $\Del_{23}$ $\Del_{12}$ \end{tabular} & $8$ & $e_{12}, e_{22}, e_{13}$ & $v_0, v_5, v_2$ \\
\hline
$t_{60}$ & \begin{tabular}{c} $\Del_{15}$ $\Del_{21}$ $\Del_{34}$ $\Del_{25}$ \\ $\Del_{16}$ $\Del_{29}$ $\Del_{37}$ $\Del_{27}$ \end{tabular} & $8$ & $e_{18}, e_{22}, e_{19}$ & $v_0, v_6, v_2$ \\
\hline
$t_{69}$ & \begin{tabular}{c} $\Del_{18}$ $\Del_{30}$ \end{tabular} & $2$ & $e_{12}, e_{32}, e_{13}$ & $v_0, v_5, v_2$ \\
\hline
$t_{80}$ & \begin{tabular}{c} $\Del_{22}$ $\Del_{33}$ $\Del_{46}$ $\Del_{35}$ \\ $\Del_{47}$ $\Del_{56}$ $\Del_{48}$ $\Del_{36}$ \end{tabular} & $8$ & $e_{33}, e_{34}, e_{19}$ & $v_7, v_6, v_2$ \\
\hline
$t_{87}$ & \begin{tabular}{c} $\Del_{26}$ $\Del_{39}$ $\Del_{49}$ $\Del_{57}$ \\ $\Del_{50}$ $\Del_{38}$ $\Del_{44}$ $\Del_{32}$ \end{tabular} & $8$ & $e_{18}, e_{37}, e_{38}$ & $v_0, v_6, v_8$ \\
\hline
$t_{93}$ & \begin{tabular}{c} $\Del_{28}$ $\Del_{40}$ $\Del_{51}$ $\Del_{41}$ \end{tabular} & $4$ & $e_{18}, e_{42}, e_{19}$ & $v_0, v_6, v_2$ \\
\hline
$t_{102}$ & \begin{tabular}{c} $\Del_{31}$ $\Del_{43}$ \end{tabular} & $2$ & $e_{44}, e_{32}, e_{45}$ & $v_0, v_9, v_2$ \\
\hline
$t_{121}$ & \begin{tabular}{c} $\Del_{42}$ $\Del_{53}$ $\Del_{58}$ $\Del_{52}$ \end{tabular} & $4$ & $e_{44}, e_{42}, e_{45}$ & $v_0, v_9, v_2$ \\
\hline
$t_{125}$ & \begin{tabular}{c} $\Del_{45}$ $\Del_{54}$ \end{tabular} & $2$ & $e_{33}, e_{50}, e_{38}$ & $v_7, v_6, v_8$ \\
\hline
$t_{133}$ & \begin{tabular}{c} $\Del_{55}$ $\Del_{59}$ \end{tabular} & $2$ & $e_{33}, e_{50}, e_{38}$ & $v_7, v_6, v_8$ \\
\hline
\end{tabular}
\caption{Triangles of type $(0\ 1\ 2)$}\label{tb:012}
\end{table}

\begin{table}[h!]
\centering
\begin{tabular}{|c|c|c|c|c|}
\hline
Triangle & Identified $\Del_j$ & Type & Edges & Vertices \\
\hline
$t_1$ & \begin{tabular}{c} $\Del_0$ $\Del_2$\end{tabular} & $(1\ 3\ 4)$ & $e_5, e_6, e_9$ & $v_1, v_3, v_4$ \\
\hline
$t_{23}$ & \begin{tabular}{c} $\Del_4$ $\Del_{10}$\end{tabular} & $(1\ 2\ 4)$ & $e_4, e_6, e_{16}$ & $v_1, v_2, v_4$ \\
\hline
$t_{28}$ & \begin{tabular}{c} $\Del_6$ $\Del_{10}$\end{tabular} & $(0\ 1\ 3)$ & $e_0, e_{17}, e_5$ & $v_0, v_1, v_3$ \\
\hline
$t_{76}$ & \begin{tabular}{c} $\Del_{22}$ $\Del_{33}$\end{tabular} & $(0\ 2\ 4)$ & $e_{34}, e_{36}, e_{16}$ & $v_7, v_2, v_4$ \\
\hline
$t_{86}$ & \begin{tabular}{c} $\Del_{26}$ $\Del_{32}$\end{tabular} & $(0\ 2\ 3)$ & $e_{37}, e_{17}, e_{39}$ & $v_0, v_8, v_3$ \\
\hline
$t_{97}$ & \begin{tabular}{c} $\Del_{31}$ $\Del_{42}$\end{tabular} & $(1\ 3\ 4)$ & $e_{46}, e_{47}, e_9$ & $v_9, v_3, v_4$ \\
\hline
$t_{126}$ & \begin{tabular}{c} $\Del_{47}$ $\Del_{56}$\end{tabular} & $(0\ 2\ 4)$ & $e_{34}, e_{36}, e_{49}$ & $v_7, v_2, v_4$ \\
\hline
$t_{129}$ & \begin{tabular}{c} $\Del_{50}$ $\Del_{57}$\end{tabular} & $(0\ 2\ 3)$ & $e_{37}, e_{48}, e_{39}$ & $v_0, v_8, v_3$ \\
\hline
$t_{131}$ & \begin{tabular}{c} $\Del_{52}$ $\Del_{58}$\end{tabular} & $(0\ 1\ 3)$ & $e_{44}, e_{48}, e_{46}$ & $v_0, v_9, v_3$ \\
\hline
$t_{132}$ & \begin{tabular}{c} $\Del_{53}$ $\Del_{58}$\end{tabular} & $(1\ 2\ 4)$ & $e_{45}, e_{47}, e_{49}$ & $v_9, v_2, v_4$ \\
\hline
\end{tabular}
\caption{Triangles with dihedral angle $\pi / 2$ and degree $2$}\label{tb:BadTri}
\end{table}

\section{The orbifold locus}\label{sec:Locus}

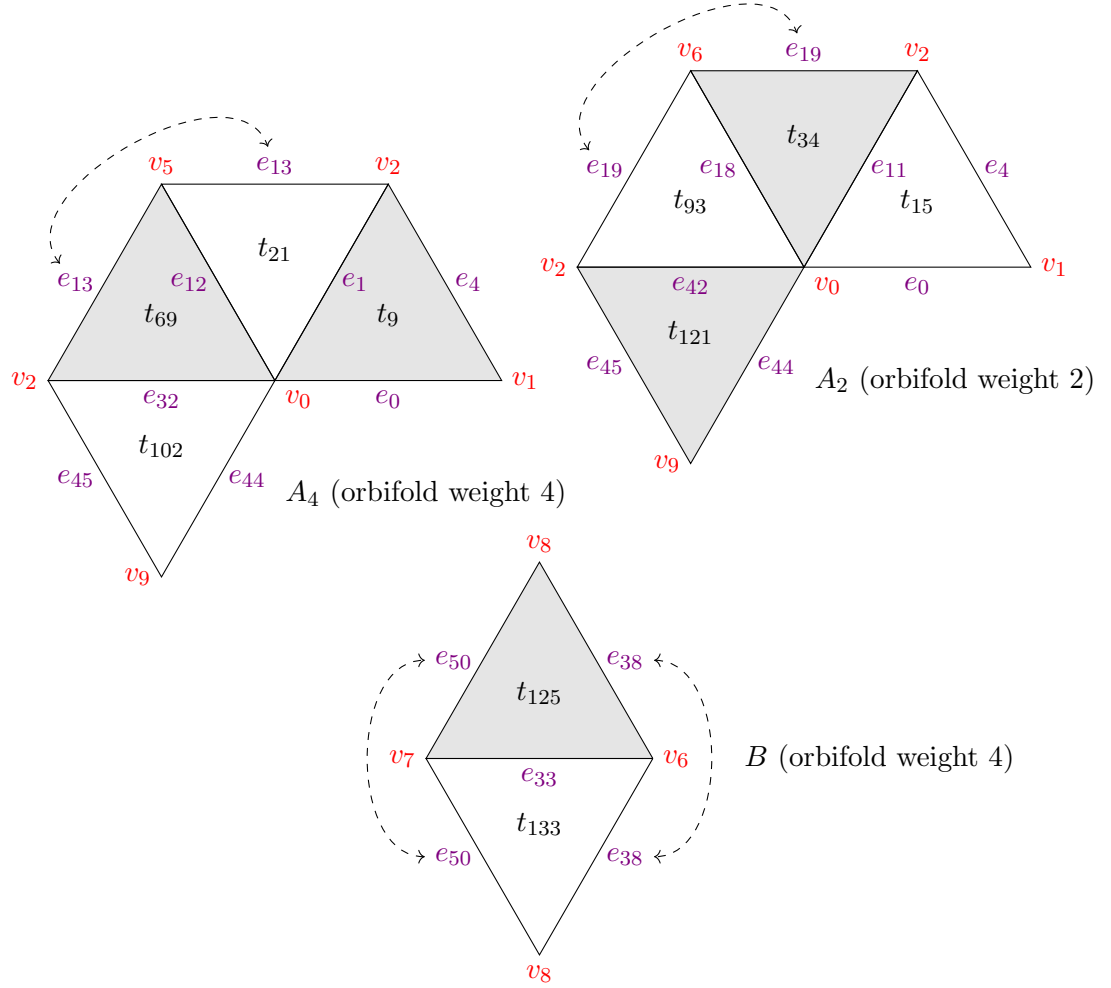
\begin{figure}[h!]
\centering
\begin{tikzpicture}
\node[red] (v0) at (0,0) {};
\node[red, below right] (v0') at (v0) {$v_0$};
\node[red] (v1) at (3,0) {};
\node[red, right] (v1') at (v1) {$v_1$};
\node[red] (v2) at (1.5, 2.59808) {};
\node[red, above] (v2') at (v2) {$v_2$};
\node[red] (v5) at (-1.5, 2.59808) {};
\node[red, above] (v5') at (v5) {$v_5$};
\node[red] (v2b) at (-3, 0) {};
\node[red, left] (v2b') at (v2b) {$v_2$};
\node[red] (v9) at (-1.5, -2.59808) {};
\node[red, left] (v9') at (v9) {$v_9$};
\draw[fill=gray, fill opacity=0.2] (0,0) -- (3,0) -- (1.5, 2.59808) -- (0,0);
\draw (0,0) -- (1.5, 2.59808) -- (-1.5, 2.59808) -- (0,0);
\draw[fill=gray, fill opacity=0.2] (0,0) -- (-3,0) -- (-1.5, 2.59808) -- (0,0);
\draw (0,0) -- (-3,0) -- (-1.5, -2.59808) -- (0,0);
\node (t9) at (barycentric cs:v0=1,v1=1,v2=1) {$t_9$};
\node (t21) at (barycentric cs:v0=1,v2=1,v5=1) {$t_{21}$};
\node (t69) at (barycentric cs:v0=1,v2b=1,v5=1) {$t_{69}$};
\node (t102) at (barycentric cs:v2b=1,v9=1,v0=1) {$t_{102}$};
\node[violet, below] (e0) at (barycentric cs:v0=1,v1=1) {$e_0$};
\node[violet, right] (e1) at (barycentric cs:v0=1,v2=1) {$e_1$};
\node[violet, right] (e4) at (barycentric cs:v1=1,v2=1) {$e_4$};
\node[violet, above] (e13) at (barycentric cs:v5=1,v2=1) {$e_{13}$};
\node[violet, left] (e12) at (barycentric cs:v0=1,v5=1) {$e_{12}$};
\node[violet, left] (e13a) at (barycentric cs:v5=1,v2b=1) {$e_{13}$};
\node[violet, below] (e32) at (barycentric cs:v0=1,v2b=1) {$e_{32}$};
\node[violet, left] (e45) at (barycentric cs:v9=1,v2b=1) {$e_{45}$};
\node[violet, right] (e44) at (barycentric cs:v9=1,v0=1) {$e_{44}$};
\draw[<->, dashed] (e13a) to [bend left = 100] (e13);
\node (weight1) at (2,-1.5) {$A_4$ (orbifold weight $4$)};
\node[red] (v0c) at (7,1.5) {};
\node[red, below right] (v0c') at (v0c) {$v_0$};
\node[red] (v1c) at (10,1.5) {};
\node[red, right] (v1c') at (v1c) {$v_1$};
\node[red] (v2c) at (8.5, 4.09808) {};
\node[red, above] (v2c') at (v2c) {$v_2$};
\node[red] (v5c) at (5.5, 4.09808) {};
\node[red, above] (v5c') at (v5c) {$v_6$};
\node[red] (v2d) at (4, 1.5) {};
\node[red, left] (v2d') at (v2d) {$v_2$};
\node[red] (v9c) at (5.5, -1.09808) {};
\node[red, left] (v9c') at (v9c) {$v_9$};
\draw (7,1.5) -- (10,1.5) -- (8.5, 4.09808) -- (7,1.5);
\draw[fill=gray, fill opacity=0.2] (7,1.5) -- (8.5, 4.09808) -- (5.5, 4.09808) -- (7,1.5);
\draw (7,1.5) -- (4,1.5) -- (5.5, 4.09808) -- (7,1.5);
\draw[fill=gray, fill opacity=0.2] (7,1.5) -- (4,1.5) -- (5.5, -1.09808) -- (7,1.5);
\node (t9c) at (barycentric cs:v0c=1,v1c=1,v2c=1) {$t_{15}$};
\node (t21c) at (barycentric cs:v0c=1,v2c=1,v5c=1) {$t_{34}$};
\node (t69c) at (barycentric cs:v0c=1,v2d=1,v5c=1) {$t_{93}$};
\node (t102c) at (barycentric cs:v2d=1,v9c=1,v0c=1) {$t_{121}$};
\node[violet, below] (e0c) at (barycentric cs:v0c=1,v1c=1) {$e_0$};
\node[violet, right] (e1c) at (barycentric cs:v0c=1,v2c=1) {$e_{11}$};
\node[violet, right] (e4c) at (barycentric cs:v1c=1,v2c=1) {$e_4$};
\node[violet, above] (e13c) at (barycentric cs:v5c=1,v2c=1) {$e_{19}$};
\node[violet, left] (e12c) at (barycentric cs:v0c=1,v5c=1) {$e_{18}$};
\node[violet, left] (e13d) at (barycentric cs:v5c=1,v2d=1) {$e_{19}$};
\node[violet, below] (e32c) at (barycentric cs:v0c=1,v2d=1) {$e_{42}$};
\node[violet, left] (e45c) at (barycentric cs:v9c=1,v2d=1) {$e_{45}$};
\node[violet, right] (e44c) at (barycentric cs:v9c=1,v0c=1) {$e_{44}$};
\draw[<->, dashed] (e13d) to [bend left = 100] (e13c);
\node (weight2) at (9,0) {$A_2$ (orbifold weight $2$)};
\node[red] (v6f) at (5,-5) {};
\node[red, right] (v6f') at (v6f) {$v_6$};
\node[red] (v8f) at (3.5, 2.59808-5) {};
\node[red, above] (v8f') at (v8f) {$v_8$};
\node[red] (v7f) at (2, -5) {};
\node[red, left] (v7f') at (v7f) {$v_7$};
\node[red] (v8g) at (3.5, -2.59808-5) {};
\node[red, below] (v8g') at (v8g) {$v_8$};
\draw[fill=gray, fill opacity=0.2] (5,-5) -- (3.5, 2.59808-5) -- (2, -5) -- (5,-5);
\draw (5,-5) -- (3.5, -2.59808-5) -- (2, -5) -- (5,-5);
\node (t125) at (barycentric cs:v6f=1,v8f=1,v7f=1) {$t_{125}$};
\node (t133) at (barycentric cs:v6f=1,v7f=1,v8g=1) {$t_{133}$};
\node[violet, right] (e38f) at (barycentric cs:v6f=1,v8f=1) {$e_{38}$};
\node[violet, left] (e50f) at (barycentric cs:v8f=1,v7f=1) {$e_{50}$};
\node[violet, below] (e33) at (barycentric cs:v7f=1,v6f=1) {$e_{33}$};
\node[violet, left] (e50g) at (barycentric cs:v7f=1,v8g=1) {$e_{50}$};
\node[violet, right] (e38g) at (barycentric cs:v6f=1,v8g=1) {$e_{38}$};
\draw[<->, dashed] (e50g) to [bend left = 90] (e50f);
\draw[<->, dashed] (e38f) to [bend left = 90] (e38g);
\node (weight3) at (8,-5) {$B$ (orbifold weight $4$)};
\end{tikzpicture}
\caption{The $(0\ 1\ 2)$ orbifold locus}\label{fig:012fig}
\end{figure}

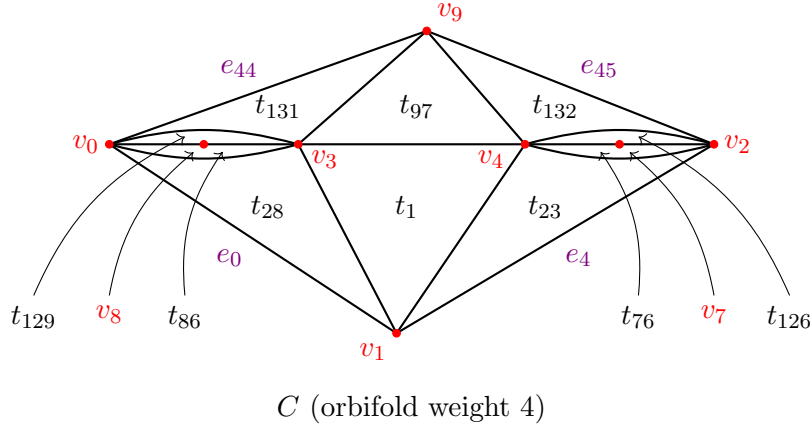
\begin{figure}[h!]
\centering
\begin{tikzpicture}
\node (v0) at (-4,0) {};
\node (v2) at (4,0) {};
\node (v1) at (-0.2,-2.5) {};
\node (v9) at (0.2,1.5) {};
\node (v4) at (1.5,0) {};
\node (v3) at (-1.5,0) {};
\node (v7) at (2.75,0) {};
\node (v8) at (-2.75,0) {};
\draw[thick] (-4,0) to (-0.2,-2.5) to (4,0);
\draw[thick] (-4,0) to (0.2,1.5) to (4,0);
\draw[thick] (-4,0) to (4,0);
\draw[thick] (-4,0) to [bend left = 15] (-1.5,0);
\draw[thick] (4,0) to [bend right = 15] (1.5,0);
\draw[thick] (4,0) to [bend left = 15] (1.5,0);
\draw[thick] (-4,0) to [bend right = 15] (-1.5,0);
\draw[thick] (-0.2,-2.5) to (1.5,0);
\draw[thick] (-0.2,-2.5) to (-1.5,0);
\draw[thick] (0.2,1.5) to (1.5,0);
\draw[thick] (0.2,1.5) to (-1.5,0);
\draw[red, fill=red] (v0) circle (0.05cm) node [left] {$v_0$};
\draw[red, fill=red] (v2) circle (0.05cm) node [right] {$v_2$};
\draw[red, fill=red] (v9) circle (0.05cm) node [above right] {$v_9$};
\draw[red, fill=red] (v1) circle (0.05cm) node [below left] {$v_1$};
\draw[red, fill=red] (v3) circle (0.05cm) node [xshift=0.35cm, yshift=-0.2cm] {$v_3$};
\draw[red, fill=red] (v4) circle (0.05cm) node [xshift=-0.45cm, yshift=-0.2cm] {$v_4$};
\draw[red, fill=red] (v7) circle (0.05cm);
\draw[red, fill=red] (v8) circle (0.05cm);
\draw[->] (4,-2) node [red, below] {$v_7$} to [bend right = 20] (v7);
\draw[->] (-4,-2) node [red, below] {$v_8$} to [bend left = 20] (v8);
\draw[->] (-3,-2) node [below] {$t_{86}$} to [bend left = 20] (-2.5,-0.1);
\draw[->] (-5,-2) node [below] {$t_{129}$} to [bend left = 20] (-3,0.1);
\draw[->] (3,-2) node [below] {$t_{76}$} to [bend right = 20] (2.5,-0.1);
\draw[->] (5,-2) node [below] {$t_{126}$} to [bend right = 20] (3,0.1);
\node[violet, below left] (e0) at (barycentric cs:v0=1,v1=1) {$e_0$};
\node[violet, below right] (e4) at (barycentric cs:v1=1,v2=1) {$e_4$};
\node[violet, above right] (e45) at (barycentric cs:v2=1,v9=1) {$e_{45}$};
\node[violet, above left] (e44) at (barycentric cs:v0=1,v9=1) {$e_{44}$};
\node (t28) at (barycentric cs:v0=1,v1=1,v3=1) {$t_{28}$};
\node (t1) at (barycentric cs:v4=1,v1=1,v3=1) {$t_1$};
\node (t23) at (barycentric cs:v4=1,v1=1,v2=1) {$t_{23}$};
\node (t132) at (barycentric cs:v4=1,v2=1,v9=1) {$t_{132}$};
\node (t97) at (barycentric cs:v4=1,v3=1,v9=1) {$t_{97}$};
\node (t131) at (barycentric cs:v0=1,v3=1,v9=1) {$t_{131}$};
\node (weight3) at (0,-3.5) {$C$ (orbifold weight $4$)};
\end{tikzpicture}
\caption{The remaining orbifold locus}\label{fig:Remaining}
\end{figure}

\begin{figure}[h!]
\centering
\begin{tikzpicture}
\node (v0) at (-4,0) {};
\node (v2) at (4,0) {};
\node (v5) at (0,5) {};
\node (v1) at (-0.2,-2.5) {};
\node (v9) at (0.2,1.5) {};
\node (v6) at (0,-6) {};
\node (v4) at (1.5,0) {};
\node (v3) at (-1.5,0) {};
\node (v7) at (2.75,0) {};
\node (v8) at (-2.75,0) {};
\draw[thick] (-0.2,-2.5) to [bend right = 20] (3.03,-2);
\draw[thick, loosely dashed] (0.2,1.5) to [bend left = 30] (3.03,-2);
\draw[thick] (-0.2,-2.5) to [bend left = 20] (-3.03,-2);
\draw[thick, loosely dashed] (0.2,1.5) to [bend right = 30] (-3.03,-2);
\draw[white, fill=gray, opacity=0.2] (4,0) to (-0.2,-2.5) to (-4,0) to (0.2,1.5) to (4,0);
\draw[white, fill=blue, opacity=0.15] (0,-6) to (-2.75,0) to [bend right = 40] (2.75,0) to (0,-6);
\draw[thick, loosely dashed] (0.2,1.5) to (0,-6);
\draw[thick] (-4,0) to (-0.2,-2.5) to (4,0);
\draw[thick, loosely dashed] (-4,0) to (0.2,1.5) to (4,0);
\draw[thick] (-4,0) to [bend right=20] (4,0);
\draw[thick, loosely dashed] (-4,0) to [bend left=58] (4,0);
\draw[thick] (-4,0) to [bend left = 10] (0,5);
\draw[thick] (4,0) to [bend right = 10] node [above right] {$A_4$} (0,5);
\draw[thick, dashed] (-4,0) to node [above left] {$C$} (4,0);
\draw[dashed, thick] (-4,0) to [bend left = 15] (-1.5,0);
\draw[dashed, thick] (4,0) to [bend right = 15] (1.5,0);
\draw[dashed, thick] (4,0) to [bend left = 15] (1.5,0);
\draw[dashed, thick] (-4,0) to [bend right = 15] (-1.5,0);
\draw[dashed, thick] (-0.2,-2.5) to (1.5,0);
\draw[dashed, thick] (-0.2,-2.5) to (-1.5,0);
\draw[dashed, thick] (0.2,1.5) to (1.5,0);
\draw[dashed, thick] (0.2,1.5) to (-1.5,0);
\draw (-4,0) to [bend right = 10] node [below left] {$A_2$} (0,-6);
\draw (4,0) to [bend left = 10] (0,-6);
\draw[thick, dashed] (0,-6) to node [below right] {$B$} (2.75,0);
\draw[thick, dashed] (0,-6) to (-2.75,0);
\draw[thick, dashed] (2.75,0) to [bend left = 40] (-2.75,0);
\draw[loosely dotted] (33/24,-3) to [bend right = 30] (-33/24,-3);
\draw[dotted] (33/24,-3) to [bend left = 40] (-33/24,-3);
\draw[thick] (-0.2,-2.5) to (0,-6);
\draw[red, fill=red] (v0) circle (0.05cm) node [left] {$v_0$};
\draw[red, fill=red] (v2) circle (0.05cm) node [right] {$v_2$};
\draw[red, fill=red] (v5) circle (0.05cm) node [above] {$v_5$};
\draw[red, fill=red] (v6) circle (0.05cm) node [below] {$v_6$};
\draw[red, fill=red] (v9) circle (0.05cm) node [above right] {$v_9$};
\draw[red, fill=red] (v1) circle (0.05cm) node [below left] {$v_1$};
\draw[red, fill=red] (v3) circle (0.05cm) node [xshift=-0.1cm, yshift=0.3cm] {$v_3$};
\draw[red, fill=red] (v4) circle (0.05cm) node [xshift=0.1cm, yshift=0.3cm] {$v_4$};
\draw[red, fill=red] (v7) circle (0.05cm);
\draw[red, fill=red] (v8) circle (0.05cm);
\draw[->] (4,-2) node [red, below] {$v_7$} to [bend right = 20] (v7);
\draw[->] (-4,-2) node [red, below] {$v_8$} to [bend left = 20] (v8);
\end{tikzpicture}
\caption{The full orbifold locus}\label{fig:FullLocus}
\end{figure}
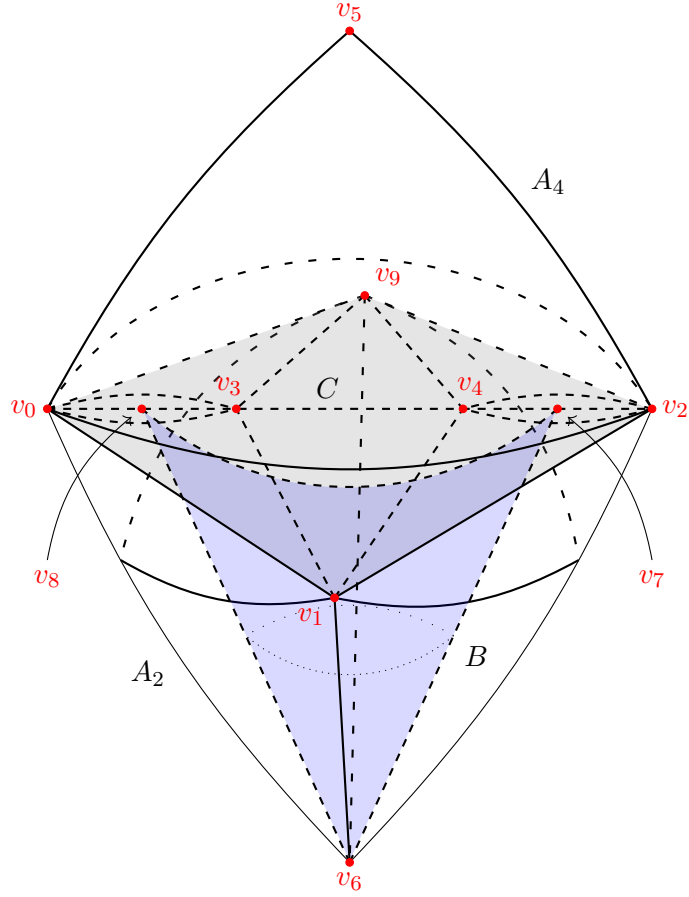

The proof of \Cref{thm:Orbifold} also gives most of the information needed to compute the orbifold locus with weights. The orbifold locus consists of four pieces that are complete orbifold quotients of $\bbH^2$. These are named $A_4$, $A_2$, $B$, and $C$, and the gluings of triangles for each are depicted in Figures \ref{fig:012fig} and \ref{fig:Remaining}, along with the orbifold weight of each piece, which is determined by the degree of the triangle and the dihedral angle of $\calS$ around that triangle. One checks using the geometry of the Coxeter simplex and the combinatorics of the gluings triangles in each piece are glued together with angle $\pi$. Note that $A_4$, $A_2$, and $C$ are polygons, meaning that the associated Fuchsian group contains orientation-reversing isometries and thus	 the induced totally geodesic suborbifold is not oriented. The complete subcomplex of $\bbP^2$ determined by the orbifold locus is shown in \Cref{fig:FullLocus}.

Tracking common vertices through the gluings, the resulting complex has ten vertices, which are recorded in \Cref{tb:Verts}, where $\Del_j(k) \cdots \Del_{j+r}(k)$ indicates that vertex $\Del_{j+\ell}(k)$ is identified with the given vertex for each $0 \le \ell \le r$. Considering $\calS$ as the unoriented quotient of $\bbH^4$ by the Coxeter group, the orbifold weight of each vertex is the order of the finite Coxeter group obtained by deleting the appropriate vertex from the Coxeter diagram. Knowing the degree of a vertex in the simplicial decomposition of $\bbP^2$ then determines the orbifold weight of that point. Orbifold weights of edges are omitted.

\begin{table}[h!]
\centering
\begin{tabular}{|c|c|c|}
\hline
Vertex & Identified $\Del_j(k)$ & Degree \\
\hline
$v_0$ & \begin{tabular}{c} $\Del_0(0) \cdots \Del_{21}(0)$ $\Del_{23}(0) \cdots \Del_{32}(0)$ \\ $\Del_{34}(0)$ $\Del_{37}(0) \cdots \Del_{44}(0)$ \\ $\Del_{49}(0) \cdots \Del_{53}(0)$ $\Del_{57}(0)$ $\Del_{58}(0)$ \end{tabular} & $48$ \\
\hline
$v_1$ & $\Del_0(1)$ $\Del_1(1)$ $\Del_2(1)$ $\Del_4(1)$ $\Del_6(1)$ $\Del_{10}(1)$ & $6$ \\
\hline
$v_2$ & \begin{tabular}{c} $\Del_0(2) \cdots \Del_{25}(2)$ $\Del_{27}(2) \cdots \Del_{31}(2)$ \\ $\Del_{33}(2) \cdots \Del_{37}(2)$ $\Del_{40}(2) \cdots \Del_{43}(2)$ \\ $\Del_{46}(2) \cdots \Del_{48}(2)$ $\Del_{51}(2) \cdots \Del_{53}(2)$ \\ $\Del_{56}(2)$ $\Del_{58}(2)$ \\ \end{tabular} & $48$ \\
\hline
$v_3$ & $\Del_0(3) \cdots \Del_{59}(3)$ & $60$ \\
\hline
$v_4$ & $\Del_0(4) \cdots \Del_{59}(4)$ & $60$ \\
\hline
$v_5$ & \begin{tabular}{c} $\Del_3(1)$ $\Del_5(1)$ $\Del_8(1)$ $\Del_9(1)$ $\Del_{12}(1)$ $\Del_{13}(1)$ \\ $\Del_{17}(1)$ $\Del_{18}(1)$ $\Del_{19}(1)$ $\Del_{23}(1)$ $\Del_{24}(1)$ $\Del_{30}(1)$ \end{tabular} & $12$ \\
\hline
$v_6$ & \begin{tabular}{c} $\Del_{7}(1)$ $\Del_{11}(1)$ $\Del_{14}(1) \cdots \Del_{16}(1)$ $\Del_{20}(1) \cdots \Del_{22}(1)$ \\ $\Del_{25}(1) \cdots \Del_{29}(1)$ $\Del_{32}(1) \cdots \Del_{41}(1)$ \\ $\Del_{44}(1) \cdots \Del_{51}(1)$ $\Del_{54}(1) \cdots \Del_{57}(1)$ $\Del_{59}(1)$ \end{tabular} & $36$ \\
\hline
$v_7$ & \begin{tabular}{c} $\Del_{22}(0)$ $\Del_{33}(0)$ $\Del_{35}(0)$ $\Del_{36}(0)$ \\ $\Del_{45}(0) \cdots \Del_{48}(0)$ $\Del_{54}(0) \cdots \Del_{56}(0)$ $\Del_{59}(0)$ \end{tabular} & $12$ \\
\hline
$v_8$ & \begin{tabular}{c} $\Del_{26}(2)$ $\Del_{32}(2)$ $\Del_{38}(2)$ $\Del_{39}(2)$ $\Del_{44}(2)$ $\Del_{45}(2)$ \\ $\Del_{49}(2)$ $\Del_{50}(2)$ $\Del_{54}(2)$ $\Del_{55}(2)$ $\Del_{57}(2)$ $\Del_{59}(2)$ \end{tabular} & $12$ \\
\hline
$v_9$ & $\Del_{31}(1)$ $\Del_{42}(1)$ $\Del_{43}(1)$ $\Del_{52}(1)$ $\Del_{53}(1)$ $\Del_{58}(1)$ & $6$ \\
\hline
\end{tabular}
\caption{Vertices in the decomposition of $\bbP^2$}\label{tb:Verts}
\end{table}

What remains is to understand how the singular locus is related to the topology of $\bbP^2$. Using linear algebra software to compute the homology of $\bbP^2$ from the given simplicial decomposition leads to the following. A file containing the various matrices needed to check the following lemma can be found at \cite{StoverCode}.

\begin{lem}\label{lem:Homology}
Let $[L]$ be a generator for $H_2(\bbP^2; \bbZ)$. Then:
\begin{itemize}

\item[$\star$] $[B] = [L]$

\item[$\star$] $[A_4 + A_2] = 2 [L]$

\item[$\star$] $[A_4 + C] = 4 [L]$

\item[$\star$] $[-A_4 + C] = 2 [L]$

\end{itemize}
In particular, $A_4 + A_2$ is homologically a smooth conic and $B$ is homologically a line tangent to the conic.
\end{lem}

It would be interesting to realize the orbifold more precisely in standard homogeneous coordinates on $\bbP^2$.

\bibliography{H4CP2}

\end{document}